\documentclass[11pt]{article}
\usepackage{wasysym,amssymb,eufrak,indentfirst,cite,bibspacing}
\begin{document}
\baselineskip 16pt
\title{On weakly S-embedded subgroups and weakly $\tau$-embedded subgroups\thanks{Research is supported by a NNSF grant of China (grant \#11071229) and Research Fund for the Doctoral Program of Higher Education of China(Grant 20113402110036).}}
\author{Xiaoyu Chen, Wenbin Guo\thanks{Corresponding author.}\\
{\small Department of Mathematics, University of Science and Technology of China,}\\ {\small Hefei 230026, P. R. China}\\
 {\small E-mail: jelly@mail.ustc.edu.cn, $\,$wbguo@ustc.edu.cn}
}
\date{}
\maketitle
\begin{abstract}
Let $G$ be a finite group. A subgroup $H$ of $G$ is said to be weakly S-embedded in $G$ if there exists $K\unlhd G$ such that $HK$ is S-quasinormal in $G$ and $H\cap K\leq H_{seG}$, where $H_{seG}$ is the subgroup generated by all those subgroups of $H$ which are S-quasinormally embedded in $G$. We say that $H$ is weakly $\tau$-embedded in $G$ if there exists $K\unlhd G$ such that $HK$ is S-quasinormal in $G$ and $H\cap K\leq H_{\tau G}$, where $H_{\tau G}$ is the subgroup generated by all those subgroups of $H$ which are $\tau$-quasinormal in $G$. In this paper, we study the properties of the weakly S-embedded subgroups and the weakly $\tau$-embedded subgroups, and use them to determine the structure of finite groups.
\end{abstract}
\renewcommand{\thefootnote}{\empty}

\section{Introduction}
Throughout this paper, all groups mentioned are finite and $G$ always denotes a finite group. All unexplained notation and terminology are standard, as in \cite{Hup,Rob,Guo4}.\par
Recall that a subgroup $H$ of $G$ is said to be \textit{S-quasinormal} (\textit{S-permutable, \textup{or} $\pi$-quasinormal}) in $G$ if $H$ permutes with every Sylow subgroup of $G$. This concept was introduced by Kegel \cite{Keg} in 1962 and has been investigated by many authors. Ballester-Bolinches and Pedraza-Aguilera \cite{Bal} extended the notion to S-quasinormally embedding in 1998: a subgroup $H$ of $G$ is said to be \textit{S-quasinormally embedded} in $G$ if every Sylow subgroup of $H$ is a Sylow subgroup of some S-quasinormal subgroup of $G$. Also, some authors considered in another way. For example, Li et al. \cite{Li} introduced the notion of SS-quasinormality in 2008: a subgroup $H$ of $G$ is said to be \textit{SS-quasinormal} in $G$ if there exists a supplement $B$ of $H$ to $G$ such that $H$ permutes with every Sylow subgroup of $B$. Chen \cite{Che} introduced the notion of S-semipermutability in 1987: a subgroup $H$ of $G$ is said to be \textit{S-semipermutable} in $G$ if $H$ permutes with every Sylow $p$-subgroup of $G$ such that $(p,|H|)=1$. Besides, V. O. Lukyanenko and A. N. Skiba \cite{Luk} introduced the notion of $\tau$-quasinormality in 2008: a subgroup $H$ of $G$ is said to be \textit{$\tau$-quasinormal} in $G$ if $H$ permutes with all Sylow $q$-subgroups $Q$ of $G$ such that $(q,|H|)=1$ and $(|H|,|Q^G|)\neq 1$.\par
It is easy to see that every SS-quasinormal subgroup and every S-semipermutable subgroup of $G$ are both $\tau$-quasinormal in $G$. In fact, it is clear that S-semipermutability implies $\tau$-quasinormality by definition. Now assume that $H$ is SS-quasinormal in $G$. Then $G$ has a subgroup $B$ such that $G=HB$ and $H$ permutes with every Sylow subgroup of $B$. Let $P$ be a Sylow $p$-subgroup of $G$ with $(p,|H|)=1$. Then there exists an element $h\in H$ such that $P^h\leq B$. Therefore $HP^h=P^hH$, and so $HP=PH$. Hence $H$ is S-semipermutable and thus $\tau$-quasinormal in $G$.\par
In 2009, Guo, Shum and A. N. Skiba \cite{Guo1} gave the definition of S-embedded subgroups: a subgroup $H$ of $G$ is said to be \textit{S-embedded} in $G$ if there exists a normal subgroup $K$ of $G$ such that $HK$ is S-quasinormal in $G$ and $H\cap K\leq H_{sG}$, where $H_{sG}$ is the subgroup generated by all those subgroups of $H$ which are S-quasinormal in $G$.\par
As a continuation of the research of S-quasinormally embedding, weakly S-embedding was introduced by Li et al. \cite{Li1} in 2011 as follows.\par
\noindent\textbf{Definition 1.1}\ \ A subgroup $H$ of $G$ is said to be \textit{weakly S-embedded} in $G$ if there exists a normal subgroup $K$ of $G$ such that $HK$ is S-quasinormal in $G$ and $H\cap K\leq H_{seG}$, where $H_{seG}$ is the subgroup generated by all those subgroups of $H$ which are S-quasinormally embedded in $G$.\par
Now we introduce the following notion:\par
\noindent\textbf{Definition 1.2}\ \ A subgroup $H$ of $G$ is said to be \textit{weakly $\tau$-embedded} in $G$ if there exists a normal subgroup $K$ of $G$ such that $HK$ is S-quasinormal in $G$ and $H\cap K\leq H_{\tau G}$, where $H_{\tau G}$ is the subgroup generated by all those subgroups of $H$ which are $\tau$-quasinormal in $G$.\par
Evidently every $\tau$-quasinormal subgroup and every S-embedded subgroup of $G$ are weakly $\tau$-embedded in $G$. Consequently, every S-quasinormal subgroup, every SS-quasinormal subgroup, and every S-semipermutable subgroup of $G$ are also weakly $\tau$-embedded in $G$. However, the next two examples show that the converse does not hold in general.\par
\noindent\textbf{Example 1.3}\ \ Let $G=S_4$ be the symmetric group of degree 4 and $H=\langle(14)\rangle$. Take $Q=\langle(123)\rangle \in Syl_3(G)$. Then clearly $Q^G=A_4$, where $A_4$ is the alternating group of degree 4. As $HQ\neq QH$, $H$ is not $\tau$-quasinormal in $G$. But since $G=HA_4$ and $H\cap A_4=1$, $H$ is weakly $\tau$-embedded in $G$.\par
\noindent\textbf{Example 1.4}\ \ Let $G=A_5$ be the alternating group of degree 5 and $H=A_4$. Because of the simplicity of $A_5$, the normal subgroups of $G$ are only $1$ and $G$. Let $P$ be a Sylow 2-subgroup of $G$ containing $\langle(15)\rangle$. As $PH\neq HP$, $H$ is not S-quasinormal in $G$, and so $H$ is not S-embedded in $G$. However, since $G=H\langle(12345)\rangle$, $H$ is $\tau$-quasinormal in $G$, and thus it is weakly $\tau$-embedded in $G$.\par
Moreover, the following examples show that weakly S-embedding and weakly $\tau$-embedding are independent of each other.\par
\noindent\textbf{Example 1.5}\ \ Let $G=A_5$ be the alternating group of degree 5, $H=\langle(123)\rangle\in Syl_3(G)$. Since every Hall $\pi$-subgroup is S-quasinormally embedded in $G$, $H$ is weakly S-embedded in $G$ (take $K=G$). On the other hand, let $Q=\langle(12345)\rangle\in Syl_5(G)$, then $Q^G=G$. It is clear that $HQ$ is not a subgroup of $G$. Hence $H$ is not $\tau$-quasinormal in $G$. In view of that $H_{\tau G}$ is $\tau$-quasinormal in $G$ (see below Lemma 2.7), we have $H$ is not weakly $\tau$-embedded in $G$.\par
\noindent\textbf{Example 1.6}\ \ Let $A=A_5$ be the alternating group of degree 5 and $B=Inn(A_5)\cong A_5$. Put $G=A\rtimes B$ and $H=A\rtimes \langle(12)(34)\rangle$. As $\pi(G)=\pi(H)$, $H$ is $\tau$-quasinormal in $G$, and thus it is weakly $\tau$-embedded in $G$. On the other hand, since $Z(A_5)=1$, $B\ntrianglelefteq G$. Clearly, the subnormal subgroups of $G$ are only $1$, $A$ and $G$. Suppose that $H$ is weakly S-embedded in $G$. Then there exists a normal subgroup $K$ of $G$ such that $HK$ is S-quasinormal in $G$ and $H\cap K\leq H_{seG}$. Since $H$ is not S-quasinormal in $G$ (see below Lemma 2.1(1)), we only can take $K=G$. This implies that $H=H_{seG}$. For any non-identity subgroup $L$ of $H$ which is S-quasinormally embedded in $G$ and any non-identity Sylow subgroup $D$ of $L$, there exists an S-quasinormal subgroup $U$ of $G$ such that $D\in Syl(U)$. In view of that $\pi(L)\subseteq \pi(|G:L|)$ and Lemma 2.1(1), we have $U=A$, and thereby $L\leq A$. It follows that $H=H_{seG}\leq A$, a contradiction. Therefore $H$ is not weakly S-embedded in $G$.\par
The purpose of this paper is to study the structure of finite groups by using the notion of weakly S-embedding and weakly $\tau$-embedding. In brief, we say a subgroup $H$ of $G$ satisfies $(\triangle)$ in $G$ if $H$ is weakly S-embedded or weakly $\tau$-embedded in $G$.\par

\section{Preliminary Lemmas}
In this paper, we use $\mathfrak{N}$, $\mathfrak{N}_{\mathfrak{p}}$ and $\mathfrak{U}$ to denote the classes of finite nilpotent, $p$-nilpotent and supersolvable groups, respectively. For a non-empty class $\mathfrak{F}$ of groups, the symbol $Z_{\mathfrak{F}}(G)$ $\bigm{(}$usually, $Z_{\mathfrak{N}}(G)$ is written as $Z_{\infty}(G)$$\bigm{)}$ denotes the $\mathfrak{F}$-hypercenter of $G$, that is, the product of all such normal subgroups $L$ of $G$ whose $G$-chief factors $H/K$ satisfy that $(H/K)\rtimes (G/C_G(H/K))\in \mathfrak{F}$. Also, the symbol $|G|_p$ denotes the order of a Sylow $p$-subgroup of $G$.\par
\noindent\textbf{Lemma 2.1} \cite{Keg,Des}. \textit{Suppose that H is S-quasinormal in G, $U\leq G$ and $N\unlhd G$. Then:}\par
\textit{\textup{(1)} H is subnormal in G.}\par
\textit{\textup{(2)} $H/H_G$ is nilpotent.}\par
\textit{\textup{(3)} $H\cap U$ is S-quasinormal in U.}\par
\textit{\textup{(4)} HN/N is S-quasinormal in G/N.}\par
\textit{\textup{(5)} If H is a p-subgroup of G, then $O^p(G)\leq N_G(H)$.}\par
\noindent\textbf{Lemma 2.2} \cite[Corollary 1]{Sch}. \textit{Suppose that A and B are S-quasinormal in G, then $\langle A,B \rangle$ and $A\cap B$ are also S-quasinormal in G.}\par
\noindent\textbf{Lemma 2.3.} \textit{Suppose that H is S-quasinormal in G, P is a Sylow p-subgroup of H, where p is a prime divisor of $|H|$. If $H_G=1$, then P is S-quasinormal in G.}\par
\noindent\textit{Proof.} This can be easily deduced from Lemma 2.1(2) and \cite[Proposition B]{Sch}.\par
\noindent\textbf{Lemma 2.4} \cite[Lemma 1]{Bal}. \textit{Suppose that H is S-quasinormally embedded in G, $U\leq G$ and $N\unlhd G$. Then:}\par
\textit{\textup{(1)} If $H\leq U$, then H is S-quasinormally embedded in U.}\par
\textit{\textup{(2)} HN is S-quasinormally embedded in G, and HN/N is S-quasinormally embedded in G/N.}\par
The following lemma can be directly obtained from Lemma 2.4.\par
\noindent\textbf{Lemma 2.5.} \textit{Suppose that $H\leq U\leq G$ and $N\unlhd G$. Then:}\par
\textit{\textup{(1)} $H_{seG}\leq H_{seU}$}.\par
\textit{\textup{(2)} $H_{seG}N/N\leq (HN/N)_{se(G/N)}$.}\par
\noindent\textbf{Lemma 2.6.} \textit{Suppose that H is $\tau$-quasinormal in G, $U\leq G$ and $N\unlhd G$.}\par
\textit{\textup{(1)} If $H\leq U$, then H is $\tau$-quasinormal in U.}\par
\textit{\textup{(2)} If $\pi(HN/N)=\pi(H)$, then HN/N is $\tau$-quasinormal in G/N.}\par
\textit{\textup{(3)} If $(|H|,|N|)=1$, then HN/N is $\tau$-quasinormal in G/N.}\par
\noindent\textit{Proof.} See \cite[Lemma 2.2]{Luk1} for (1) and (3). Now we prove (2). Let $Q/N\in Syl_q(G/N)$ such that $(|HN/N|,|(Q/N)^{(G/N)}|)\neq 1$ and $q\notin \pi(HN/N)=\pi(H)$. Then for some Sylow $q$-subgroup $G_q$ of $G$, we have $Q=G_qN$. Since $(Q/N)^{(G/N)}=Q^G/N=(G_qN)^G/N=G_q^GN/N\cong G_q^G/G_q^G\cap N$, $(|H|,|G_q^G|)\neq 1$. Hence by the hypothesis, $(HN/N)(Q/N)=HG_qN/N=G_qHN/N=(Q/N)(HN/N)$.\par
From Lemma 2.6 we directly have:\par
\noindent\textbf{Lemma 2.7.} \textit{Suppose that $H\leq U\leq G$ and $N\unlhd G$. Then:}\par
\textit{\textup{(1)} If H is a p-subgroup, then $H_{\tau G}$ is $\tau$-quasinormal in G and $H_G\leq H_{\tau G}$} \cite[Lemma 2.3(1)]{Luk1}.\par
\textit{\textup{(2)} $H_{\tau G}\leq H_{\tau U}$} \cite[Lemma 2.3(2)]{Luk1}.\par
\textit{\textup{(3)} If H is a p-subgroup, then $H_{\tau G}N/N\leq (HN/N)_{\tau (G/N)}$.}\par
\textit{\textup{(4)} If $(|H|,|N|)=1$, then $H_{\tau G}N/N\leq (HN/N)_{\tau (G/N)}$.}\par
\noindent\textbf{Lemma 2.8.} \textit{Let P be a p-subgroup of G. Then the following statements are equivalent:}\par
\textit{\textup{(1)} P is S-quasinormal in G.}\par
\textit{\textup{(2)} $P\leq O_p(G)$ and P is S-quasinormally embedded in G} \cite[Lemma 2.4]{Li4}.\par
\textit{\textup{(3)} $P\leq O_p(G)$ and P is $\tau$-quasinormal in G} \cite[Lemma 2.2(4)]{Luk1}.\par
\noindent\textbf{Lemma 2.9} \cite[Lemma 2.4]{Li1}. \textit{Suppose that H is weakly S-embedded in G, $U\leq G$ and $N\unlhd G$.}\par
\textit{\textup{(1)} If $H\leq U$, then H is weakly S-embedded in U.}\par
\textit{\textup{(2)} If $N\leq H$, then H/N is weakly S-embedded in G/N.}\par
\textit{\textup{(3)} If $(|H|,|N|)=1$, then HN/N is weakly S-embedded in G/N.}\par
Now we give some basic properties of weakly $\tau$-embedded subgroups.\par
\noindent\textbf{Lemma 2.10.} \textit{Suppose that H is weakly $\tau$-embedded in G, $U\leq G$ and $N\unlhd G$.}\par
\textit{\textup{(1)} If $H\leq U$, then H is weakly $\tau$-embedded in U.}\par
\textit{\textup{(2)} If H is a p-subgroup and $N\leq H$, then H/N is weakly $\tau$-embedded in G/N.}\par
\textit{\textup{(3)} If $(|H|,|N|)=1$, then HN/N is weakly $\tau$-embedded in G/N.}\par
\noindent\textit{Proof.} By the hypothesis, there exists a normal subgroup $K$ of $G$ such that $HK$ is $S$-quasinormal in $G$ and $H\cap K\leq H_{\tau G}$. Then:\par
(1) $K\cap U\unlhd U$ and $H(K\cap U)=HK\cap U$ is $S$-quasinormal in $U$ by Lemma 2.1(3). In view of Lemma 2.7(2), we have $H\cap (K\cap U)=H\cap K\leq H_{\tau G}\leq H_{\tau U}$. Hence $H$ is weakly $\tau$-embedded in $U$.\par
(2) $KN/N\unlhd G/N$ and $(H/N)(KN/N)=HKN/N$ is $S$-quasinormal in $G/N$ by Lemma 2.1(4). Since $(H/N)\cap (KN/N)=(H\cap K)N/N$, $(H/N)\cap (KN/N)\leq H_{\tau G}N/N\leq (H/N)_{\tau (G/N)}$ by Lemma 2.7(3). Therefore $H/N$ is weakly $\tau$-embedded in $G/N$.\par
(3) $KN/N\unlhd G/N$ and $(HN/N)(KN/N)=HKN/N$ is $S$-quasinormal in $G/N$ by Lemma 2.1(4). It is easy to see that $(|N\cap HK:N\cap H|,|N\cap HK:N\cap K|)=1$. Then $N\cap HK=(N\cap H)(N\cap K)$, and so $HN\cap KN=(H\cap K)N$ by \cite[Chapter A, Lemma (1.2)]{Doe}. It follows from Lemma 2.7(4) that $(HN/N)\cap (KN/N)=(H\cap K)N/N\leq H_{\tau G}N/N\leq (HN/N)_{\tau (G/N)}$. Consequently, $HN/N$ is weakly $\tau$-embedded in $G/N$.\par
\noindent\textbf{Lemma 2.11} \cite[Lemma 3.4.7]{Guo4}. \textit{Let $\mathfrak{F}$ be a saturated formation. If $G\not\in \mathfrak{F}$, but every proper subgroup of G belongs to $\mathfrak{F}$. Furthermore, $G$ has a normal Sylow p-subgroup $G_p\neq 1$, where p is a prime divisor of $|G|$. Then:}\par
\textit{\textup{(1)} $G_p=G^{\mathfrak{F}}$.}\par
\textit{\textup{(2)} $F(G)=F_p(G)=G_p\Phi(G)$.}\par
\noindent\textbf{Lemma 2.12.} \textit{Let p be a prime divisor of $|G|$ with $(|G|,(p-1)(p^2-1)\cdots (p^n-1))=1$. If $H\unlhd G$ with $p^{n+1}\;{\not|}\;|H|$ for some integer $n\geq 1$ and G/H is p-nilpotent, then G is p-nilpotent. In particular, if $p^{n+1}\;{\not|}\;|G|$, then G is p-nilpotent.}\par
\noindent\textit{Proof.} Assume that the result is false and let $(G,H)$ be a counterexample for which $|G|+|H|$ is minimal. For every non-trivial subgroup $F$ of $G$, $(F,H\cap F)$ satisfies the hypothesis. Then $F$ is $p$-nilpotent due to the choice of $(G,H)$. It induces that $G$ is a minimal non-$p$-nilpotent group, and thus a minimal nonnilpotent group by \cite[Chapter IV, Theorem 5.4]{Hup}. By \cite[Theorem 3.4.11]{Guo4} and Lemma 2.11, $G=P\rtimes Q$ with $P=G^{\mathfrak{N}}=G^{\mathfrak{N}_\mathfrak{p}}\in Syl_p(G)$. If $N$ is a minimal normal subgroup of $G$, then $(G/N,HN/N)$ satisfies the hypothesis. Hence $G/N$ is $p$-nilpotent. This implies that $P=G^{\mathfrak{N}_\mathfrak{p}}=N$ is an elementary abelian group. Clearly $P\leq H$, and so $|P|\leq p^n$. Since $N_G(P)/C_G(P)\apprle Aut(P)$, $|N_G(P)/C_G(P)|\Bigm|(p-1)(p^2-1)\cdots (p^n-1)$. Consequently, $N_G(P)=C_G(P)$. So by Burnside's Theorem, $G$ is $p$-nilpotent, a contradiction.\par
\noindent\textbf{Lemma 2.13} \cite[Chapter VI, Lemma 4.10]{Hup}. \textit{Let A and B be subgroups of G such that $G\neq AB$. If $A^gB=BA^g$ for all $g\in G$, then there exists a non-trivial normal subgroup N of G containing either A or B.}\par
\noindent\textbf{Lemma 2.14.} \textit{Let p be a prime divisor of $|G|$ with $(|G|,p-1)=1$. If G has a Hall $p'$-subgroup, then any two Hall $p'$-subgroups of G are conjugate in G.}\par
\noindent\textit{Proof.} If $p=2$, then by \cite[Theorem A]{Cro}, any two Hall $p'$-subgroups are conjugate in $G$. If $p$ is an odd prime, then $2{\not|}|G|$. By Feit-Thompson's Theorem, $G$ is soluble. Hence any two Hall $p'$-subgroups are conjugate in $G$.\par
The next lemma is evident.\par
\noindent\textbf{Lemma 2.15.} \textit{Let p be a prime divisor of $|G|$, $H\leq G$ and $N\unlhd G$. If $|HN/N|_p\geq p^{n+1}$ for some integer $n\geq 1$, then for every $T/N\in Syl_p(HN/N)$, there exists a Sylow p-subgroup P of H such that $T/N=PN/N$; for every n-maximal subgroup $T_n/N$ of $T/N$, there exists an n-maximal subgroup $P_n$ of $P$ such that $T_n/N=P_nN/N$ and $P_n\cap N=P\cap N$.}\par
\noindent\textbf{Lemma 2.16} \cite[Lemma 2.8]{Wan1}. \textit{Let M be a maximal subgroup of a group G, P a normal p-subgroup of G such that $G=PM$, where p is a prime divisor of $|G|$. Then $P\cap M$ is normal in G.}\par
\noindent\textbf{Lemma 2.17} \cite{Tat}. \textit{Let P be a Sylow p-subgroup of G, where p is a prime divisor of $|G|$. If $N\unlhd G$ and $N\cap P\subseteq \Phi(P)$, then N is p-nilpotent.}\par
\noindent\textbf{Lemma 2.18} \cite[Chapter A, (4.13) Proposition(b)]{Doe}. \textit{Let $G=G_1\times \cdots \times G_r$ with each $G_i$ a non-abelian simple group. Then a subgroup S is subnormal in G if and only if S is a (direct) product of a subset of the factors $G_i$.}\par
\noindent\textbf{Lemma 2.19} \cite[Lemma 2.1]{Guo6}\textbf{.} \textit{Let $\mathfrak{F}$ be a non-empty saturated formation. Suppose that $A\leq G$.}\par
\textit{\textup{(1)} If $A\unlhd G$, then $AZ_{\mathfrak{F}}(G)/A\leq Z_{\mathfrak{F}}(G/A)$.}\par
\textit{\textup{(2)} If $\mathfrak{F}$ is S-closed, then $Z_{\mathfrak{F}}(G)\cap A\leq Z_{\mathfrak{F}}(A)$.}\par
\textit{\textup{(3)} If $G\in \mathfrak{F}$, then $Z_{\mathfrak{F}}(G)=G$.}\par
\noindent\textbf{Lemma 2.20} \cite[Chapter III, Theorem 7.2]{Hup}. \textit{Let G be a p-group. Suppose that $H_1$ and $H_2$ are \textup{(}normal\textup{)} subgroups of G such that $H_1\leq H_2$ and $|H_2:H_1|=p^s$. Then for any integer $1\leq t\leq s$, there exists a \textup{(}normal\textup{)} subgroup $H_3$ of G, which satisfies that $H_1\leq H_3\leq H_2$ and $|H_3:H_1|=p^t$.}\par
\noindent\textbf{Lemma 2.21} \cite[Lemma 2.8]{Li4}. \textit{Let p be a prime dividing the order of G, G is $A_4$-free and $(|G|,p-1)=1$. Assume that N is a normal subgroup of G such that G/N is p-nilpotent and the order of N is not divisible by $p^3$. Then G is p-nilpotent.}\par
\noindent\textbf{Lemma 2.22.} \textit{Let p be a prime divisor of $|G|$ with $(|G|,p-1)=1$. If G has cyclic Sylow p-subgroups, then G is p-nilpotent.}\par
\noindent\textit{Proof.} See the proof of \cite[(10.1.9)]{Rob}.\par
\noindent\textbf{Lemma 2.23} \cite[Lemma 2.9]{Guo2}. \textit{Let $\mathfrak{F}$ be a saturated formation containing $\mathfrak{U}$ \textup{(}containing $\mathfrak{N}$\textup{)} and G a group with a normal subgroup E such that $G/E\in \mathfrak{F}$. If E is cyclic \textup{(}if E is contained in Z(G), respectively\textup{)}, then $G\in \mathfrak{F}$.}\par
\noindent\textbf{Lemma 2.24} \cite[Main Theorem]{Li1}. \textit{Let $\mathfrak{F}$ be a saturated formation containing $\mathfrak{U}$. Then $G\in \mathfrak{F}$ if and only if G has a normal subgroup E such that $G/E\in \mathfrak{F}$ and for any non-cyclic Sylow subgroup P of the generalized Fitting subgroup $F^*(E)$, every maximal subgroup of P not having a supersolvable supplement in G or every cyclic subgroup H of P of prime order or order \textup{4} $\bigm{(}$when P is a non-abelian \textup{2}-group and $H\nsubseteq  Z_{\infty}(G)$$\bigm{)}$ without a supersolvable supplement in G is weakly S-embedded in G.}\par

\section{Main Results}
Note that a subgroup $M_n$ of $G$ is said to be an $n$-maximal $(n\geq 1)$ subgroup of $G$ if $G$ has a subgroup chain: $M_{n}<M_{n-1}< \cdots <M_1<M_0=G$ such that $M_i$ is a maximal subgroup of $M_{i-1}$ $(1\leq i \leq n)$.\par
\medskip
\noindent\textbf{Theorem 3.1.} \textit{Let G be a group and p a prime divisor of $|G|$ with $(|G|,(p-1)(p^2-1)\cdots (p^n-1))=1$ for some integer $n\geq 1$. Then G is p-nilpotent if and only if there exists a normal subgroup H of G such that G/H is p-nilpotent and for any Sylow p-subgroup P of H, every n-maximal subgroup of P not containing $P\cap G^{\mathfrak{N}_\mathfrak{p}}$ \textup{(}if exists\textup{)} either has a p-nilpotent supplement in G or satisfies $(\triangle)$ in $G$.}\par
\medskip                                                                                                                                                                                                                                                                                                          \noindent\textit{Proof.} The necessity is obvious. So we need only to prove the sufficiency. Assume that the result is false and let $(G,H)$ be a counterexample for which $|G|+|H|$ is minimal. We proceed the proof via the following steps:\par
(1) \textit{${|H|}_p\geq p^{n+1}$.}\par
By Lemma 2.12, if ${|H|}_p\leq p^n$, then $G$ is $p$-nilpotent, a contradiction.\par
(2) \textit{$G$ is not a non-abelian simple group.}\par
If $G$ is a non-abelian simple group, then $G^{\mathfrak{N}_\mathfrak{p}}=H=G$ for $G^{\mathfrak{N}_\mathfrak{p}}\neq 1$. Since ${|G|}_p\geq p^{n+1}$, we may let $P_n$ be an $n$-maximal subgroup of a Sylow $p$-subgroup $P$ of $G$. By the hypothesis, $P_n$ either satisfies $(\triangle)$ or has a $p$-nilpotent supplement in $G$.\par
(i) \textit{Case \textup{1:} $P_n$ satisfies $(\triangle)$ in $G$.}\par
In this case, there exists a normal subgroup $K$ of $G$ such that $P_nK$ is S-quasinormal in $G$ and $P_n\cap K\leq (P_n)_{seG}$ or $P_n\cap K\leq (P_n)_{\tau G}$. As $G$ is simple, $K=1$ or $G$.\par
Suppose that $K=1$. Then $P_n$ is S-quasinormal and so subnormal in $G$ by Lemma 2.1(1). Thus $P_n=1$ for $G\neq P_n$. It follows that ${|G|}_p=p^n$, which contradicts (1). Now assume that $K=G$. Then $P_n=(P_n)_{seG}$ or $P_n=(P_n)_{\tau G}$. In the former case, for any subgroup $D$ of $P_n$ which is S-quasinormally embedded in $G$, there exists an S-quasinormal subgroup $U$ of $G$ such that $D\in Syl_p(U)$. Since $G\neq U$, $D=U=1$. Therefore $P_n=(P_n)_{seG}=1$, and so ${|G|}_p=p^n$, a contradiction. In the latter case, $P_n$ is $\tau$-quasinormal in $G$ by Lemma 2.7(1). As $G$ is not a $p$-group, we can take a prime $q\neq p$ dividing $|G|$ and a Sylow $q$-subgroup $Q$ of $G$. Suppose that $P_n\neq 1$. As $Q^G=G$, $(p,|Q^G|)\neq 1$. So we have $Q^gP_n=P_nQ^g$ for all $g\in G$. If $G=P_nQ$, then $G$ is soluble by Burnside's $p^aq^b$-Theorem, a contradiction. Hence $G\neq P_nQ$. By Lemma 2.13, there exists a non-trivial normal subgroup $X$ of $G$ such that either $P_n\leq X$ or $Q\leq X$, which is impossible. So we obtain that $P_n=1$ and ${|G|}_p=p^n$, a contradiction again.\par
(ii) \textit{Case \textup{2:} $P_n$ has a $p$-nilpotent supplement $T$ in $G$.}\par
Let $T_{p'}$ be the normal $p$-complement of $T$. Then $G=P_nT=P_nN_G(T_{p'})$. If $N_G(T_{p'})=G$, then $T_{p'}\unlhd G$, and so $T_{p'}=1$ or $G$. This implies that $G$ is $p$-nilpotent, a contradiction. Thus $N_G(T_{p'})< G$. Obviously $P\cap N_G(T_{p'})\in Syl_p(N_G(T_{p'}))$ and $P\cap N_G(T_{p'})<P$. Let $P_2$ be a maximal subgroup of $P$ containing $P\cap N_G(T_{p'})$ and $P_{n2}$ an $n$-maximal subgroup of $P$ contained in $P_2$. If $P_{n2}$ satisfies $(\triangle)$ in $G$, the same discussion as above shows that it is impossible. Hence $P_{n2}$ has a $p$-nilpotent supplement $E$ in $G$. Let $E_{p'}$ be the normal $p$-complement of $E$, then we also have that $G=P_{n2}N_G(E_{p'})=P_2N_G(E_{p'})$. Since $T_{p'}$ and $E_{p'}$ are Hall $p'$-subgroups of $G$, there exists an element $g\in P_2$ such that $T_{p'}={(E_{p'})}^g$ by Lemma 2.14. Consequently, $G=(P_2N_G(E_{p'}))^g=P_2N_G(T_{p'})$, and thereby $P=P_2(P\cap N_G(T_{p'}))=P_2$, a contradiction. This completes the proof of (2).\par
(3) \textit{If $1\neq L\unlhd G$ such that $L\leq H$ or $(|L|,p)=1$, then $G/L$ is $p$-nilpotent.}\par
If ${|HL/L|}_p\leq p^n$, then $G/L$ is $p$-nilpotent owing to Lemma 2.12. So we may assume that ${|HL/L|}_p\geq p^{n+1}$. Let $PL/L$ be a Sylow $p$-subgroup of $HL/L$ with $P\in Syl_p(H)$ and $P_nL/L$ an $n$-maximal subgroup of $PL/L$ not containing $(PL/L)\cap (G/L)^{\mathfrak{N}_\mathfrak{p}}$ such that $P_n$ is an $n$-maximal subgroup of $P$ and $P_n\cap L=P\cap L$ (see Lemma 2.15). As $G^{\mathfrak{N}_\mathfrak{p}}\leq H$, $(|PG^{\mathfrak{N}_\mathfrak{p}}\cap L:P\cap L|,|PG^{\mathfrak{N}_\mathfrak{p}}\cap L:G^{\mathfrak{N}_\mathfrak{p}}\cap L|)=1$. Thus $PL\cap G^{\mathfrak{N}_\mathfrak{p}}L=(P\cap G^{\mathfrak{N}_\mathfrak{p}})L$ by \cite[Chapter A, Lemma (1.2)]{Doe}. If $P\cap G^{\mathfrak{N}_\mathfrak{p}}\leq P_n$, then $(PL/L)\cap (G/L)^{\mathfrak{N}_\mathfrak{p}}=(PL/L)\cap (G^{\mathfrak{N}_\mathfrak{p}}L/L)=(P\cap G^{\mathfrak{N}_\mathfrak{p}})L\leq P_nL/L$, a contradiction. So by the hypothesis, $P_n$ either has a $p$-nilpotent supplement $T$ or satisfies $(\triangle)$ in $G$. In the former case, $P_nL/L$ has a $p$-nilpotent supplement $TL/L$ in $G/L$. In the latter case, there exists a normal subgroup $K$ of $G$ such that $P_nK$ is S-quasinormal in $G$ and $P_n\cap K\leq (P_n)_{seG}$ or $P_n\cap K\leq (P_n)_{\tau G}$. Note that $|L\cap P_nK:L\cap K|=|K(L\cap P_nK):K|\Bigm||P_nK:K|$ and $|L\cap P_nK:L\cap P_n|\Bigm||L:L\cap P_n|=|L:L\cap P|=|PL:P|$. Since $L\leq H$ or $(|L|,p)=1$, $p\;{\not|}\;|PL:P|$. This implies that $(|L\cap P_nK:L\cap K|,|L\cap P_nK:L\cap P_n|)=1$. Hence $P_nL\cap KL=(P_n\cap K)L$ as above. In view of Lemma 2.5(2) and Lemma 2.7(3), either $(P_nL/L)\cap (KL/L)\leq (P_n)_{seG}L/L\leq (P_nL/L)_{seG}$ or $(P_nL/L)\cap (KL/L)\leq (P_n)_{\tau G}L/L\leq (P_nL/L)_{\tau G}$. This shows that $P_nL/L$ satisfies $(\triangle)$ in $G/L$. Therefore $(G/L,HL/L)$ satisfies the hypothesis. By the choice of $(G,H)$, $G/L$ is $p$-nilpotent.\par
(4) \textit{If $P\leq F<G$, then $F$ is $p$-nilpotent. In particular, if $H<G$, then $H$ is $p$-nilpotent.}\par
Obviously $P\in Syl_p(H\cap F)$. Let $P_n$ be an $n$-maximal subgroup of $P$ not containing $P\cap F^{\mathfrak{N}_\mathfrak{p}}$. Suppose that $P\cap G^{\mathfrak{N}_\mathfrak{p}}\leq P_n$. Then $P\cap F^{\mathfrak{N}_\mathfrak{p}}\leq P\cap (G^{\mathfrak{N}_\mathfrak{p}}\cap F)\leq P_n$, which is impossible. Hence by the hypothesis, $P_n$ either has a $p$-nilpotent supplement $T$ or satisfies $(\triangle)$ in $G$. It induces from Lemma 2.9(1) and Lemma 2.10(1) that $P_n$ either has a $p$-nilpotent supplement $T\cap F$ or satisfies $(\triangle)$ in $F$. In view of the conjugacy of the Sylow $p$-subgroups, $(F,H\cap F)$ satisfies the hypothesis, and so $F$ is $p$-nilpotent by the choice of $(G,H)$.\par
(5) \textit{$G$ has a unique minimal normal subgroup $N=G^{\mathfrak{N}_\mathfrak{p}}$ contained in $H$ such that $N\nleq \Phi(G)$ and ${|N|}_p\geq p^{n+1}$.}\par
It follows directly from (3) and Lemma 2.12.\par
(6) \textit{$O_{p'}(G)=1$.}\par
If not, by (3), $G/O_{p'}(G)$ is $p$-nilpotent. Consequently, $G$ is $p$-nilpotent, a contradiction.\par
(7) \textit{$O_p(H)=1$.}\par
If not, then:\par
(i) \textit{$N=O_p(H)\leq O_p(G)$, and $G=N\rtimes M$ for some maximal subgroup $M$ of $G$.}\par
Since $O_p(H)\;char\;H\unlhd G$, $N\leq O_p(H)$, and thereby $N$ is abelian. In view of $N\nleq \Phi(G)$, there exists a maximal subgroup $M$ of $G$ such that $G=N\rtimes M=O_p(H)M$. It follows from Lemma 2.16 that $O_p(H)\cap M\unlhd G$. The uniqueness of $N$ yields $O_p(H)\cap M=1$, which implies that  $O_p(H)=N(O_p(H)\cap M)=N\leq O_p(G)$.\par
(ii) \textit{Let $G_p$ be a Sylow $p$-subgroup of $G$ which satisfies that $P=G_p\cap H$. Then there exist a maximal subgroup $P_1$ of $P$ and an $n$-maximal subgroup $P_n$ of $P$ contained in $P_1$ such that $P_1\unlhd G_p$ and $P=NP_n=NP_1$. Furthermore, $N\nleq P_1$ and $N\cap P_n$ is an $n$-maximal subgroup of $N$.}\par
Clearly $G_p=N(G_p\cap M)$. As $|G_p:G_p\cap M|=|G_pM:M|=|G:M|=|N|\geq p^{n+1}$, there exist an $n$-maximal subgroup ${(G_p)}_n$ and a maximal subgroup ${(G_p)}_1$ of $G_p$ such that $G_p\cap M\leq {(G_p)}_n\leq {(G_p)}_1$. This implies that $G_p=N{(G_p)}_n=N{(G_p)}_1$. Let $P_n={(G_p)}_n\cap H$ and $P_1={(G_p)}_1\cap H$. Obviously $P_1\unlhd G_p$ and $P\cap M\leq P_n\leq P_1$. Thus $P=N(P\cap M)=NP_n=NP_1$. Moreover, since $|P:P_n|=|G_p\cap H:{(G_p)}_n\cap H|=|G_p:{(G_p)}_n|=p^n$ and $|P:P_1|=|G_p\cap H:{(G_p)}_1\cap H|=|G_p:{(G_p)}_1|=p$, $P_n$ is an $n$-maximal subgroup of $P$ and $P_1$ is a maximal subgroup of $P$. It is evident that $N\nleq P_1$. As $|N:N\cap P_n|=|P:P_n|=p^n$, $N\cap P_n$ is an $n$-maximal subgroup of $N$.\par
(iii) \textit{$P_n$ satisfies $(\triangle)$ in $G$, that is, there exists a normal subgroup $K$ of $G$ such that $P_nK$ is S-quasinormal in $G$ and $P_n\cap K\leq (P_n)_{seG}$ or $P_n\cap K\leq (P_n)_{\tau G}$.}\par
Assume that $P_n$ does not satisfy $(\triangle)$ in $G$, then $P_n$ has a $p$-nilpotent supplement $T$ in $G$ by the hypothesis. Let $T_{p'}$ be the normal $p$-complement of $T$, then $G=P_nT=P_nN_G(T_{p'})$. Since $M\cong G/N$ is $p$-nilpotent, $M$ has a normal $p$-complement $M_{p'}$ such that $M\leq N_G(M_{p'})\leq G$. If $N_G(M_{p'})=G$, then $M_{p'}\unlhd G$, and thus $M_{p'}=1$ by (6). This shows that $G$ is a $p$-group, a contradiction. Therefore $N_G(M_{p'})=M$. By Lemma 2.14, there exists an element $g\in P_n$ such that $M_{p'}={(T_{p'})}^g$. Then $G=(P_nN_G(T_{p'}))^g=P_nN_G(M_{p'})=P_nM$, and thereby $P=P_n(P\cap M)=P_n$, a contradiction. Hence (iii) holds.\par
(iv) \textit{$P_n\cap K=1$.}\par
First suppose that $P_n\cap K\leq (P_n)_{seG}$. Let $D_1,D_2,\cdots D_s$ be all subgroups of $P_n$ which are S-quasinormally embedded in $G$, and so $(P_n)_{seG}=\langle D_1,D_2,\cdots D_s \rangle$. By definition, there exist S-quasinormal subgroups $U_1,U_2,\cdots U_s$ of $G$ with $D_i\in Syl_p(U_i)$ $(1\!\leq\! i \!\leq \! s)$. Assume that $(U_i)_G\neq 1$. Let $N_1$ be a minimal normal subgroup of $G$ contained in $(U_i)_G$. If $N_1\nleq H$, then $N_1\cap H=1$. Let $N_2/H\unlhd G/H$ contained in $HN_1/H$, then $N_2\unlhd G$, and therefore $N_1\cap N_2=1$ or $N_1$. Thus $N_2=H(N_1\cap N_2)=H$ or $HN_1$. So $HN_1/H$ is a minimal normal subgroup of $G/H$. Since $G/H$ is $p$-nilpotent, it is also $p$-soluble. It follows from (6) that $N_1\cong HN_1/H$ is a $p$-subgroup . Consequently, $N_1\leq D_i\leq H$, a contradiction. Hence $N_1\leq H$, then by (5), $N_1=N\leq U_i\leq P_n$, also a contradiction.\par
Thus $(U_i)_G=1$ $(1\!\leq\! i \!\leq \! s)$. By Lemma 2.2 and Lemma 2.3, $D_i$ is S-quasinormal in $G$ and so is $(P_n)_{seG}$. Therefore, $(P_n)_{seG}\leq O_p(H)=N$ and $O^p(G)\leq N_G((P_n)_{seG})$ by Lemma 2.1(1) and Lemma 2.1(5). If $(P_n)_{seG}\neq 1$, then $N\leq ((P_n)_{seG})^G=((P_n)_{seG})^{O^p(G)G_p}=((P_n)_{seG})^{G_p}\leq (P_n\cap N)^{G_p}\leq (P_1\cap N)^{G_p}=P_1\cap N\leq N$. This implies that $N=P_1\cap N$ and $N\leq P_1$, a contradiction. Thus $(P_n)_{seG}=1$, and so $P_n\cap K=1$.\par
Now consider that $P_n\cap K\leq (P_n)_{\tau G}$. If $H\cap K\neq 1$, then $N\leq H\cap K\leq K$. Thereupon $N\cap P_n\leq P_n\cap K\leq (P_n)_{\tau G}$, and so $N\cap P_n=N\cap (P_n)_{\tau G}$. Take an arbitrary prime $q\neq p$ dividing $|G|$ and a Sylow $q$-subgroup $Q$ of $G$. Suppose that $(p,|Q^G|)=1$, then $Q^G\leq O_{p'}(G)$, which contradicts (6). By Lemma 2.7(1), $(P_n)_{\tau G}$ is $\tau$-quasinormal in $G$, and so $(P_n)_{\tau G}Q$ is a subgroup of $G$. Note that $N\cap P_n$ is an $n$-maximal subgroup of $N$. Since $|N:N\cap (P_n)_{\tau G}Q|=|N(P_n)_{\tau G}Q:(P_n)_{\tau G}Q|=|N(P_n)_{\tau G}:(P_n)_{\tau G}|=|N:N\cap P_n|=p^n$, $N\cap (P_n)_{\tau G}Q=N\cap P_n$. It follows that $Q\leq N_G(N\cap (P_n)_{\tau G}Q)=N_G(N\cap P_n)$, and thus $O^p(G)\leq N_G(N\cap P_n)$. If $N\cap P_n\neq 1$, then $N\leq (N\cap P_n)^G=(N\cap P_n)^{G_p}\leq (N\cap P_1)^{G_p}=N\cap P_1$. This implies that $N\leq P_1$, a contradiction. Hence $N\cap P_n=1$. Then $|N|=p^n$, also a contradiction. Therefore, $H\cap K=1$, and so $P_n\cap K=1$. The proof of (iv) is completed.\par
(v) \textit{$H=G$.}\par
If $H<G$, then $H$ is $p$-nilpotent by (4). Let $H_{p'}$ be the normal $p$-complement of $H$. Then $H_{p'}\unlhd G$. Hence $H_{p'}=1$ by (6), which implies that $H=N=P$. Since $P_n(H\cap K)=P_nK\cap H$ is S-quasinormal in $G$ by Lemma 2.2, we have $O^p(G)\leq N_G(P_n(H\cap K))$. If $P_n(H\cap K)\neq 1$, then $H=(P_n(H\cap K))^G=(P_n(H\cap K))^{G_p}\leq (P_1(H\cap K))^{G_p}=P_1(H\cap K)$. Therefore, $H=P_1(H\cap K)$ and $H\cap K\neq 1$. It induces that $H\cap K=H$ by the minimality of $H$. Then $H\leq K$, and so $P_n=P_n\cap K=1$ by (iv). Consequently, $|H|=p^n$, which contradicts (1). Hence $P_n(H\cap K)=1$. We also obtain that $P_n=1$, a contradiction as above.\par
(vi) \textit{Final contradiction of \textup{(7)}.}\par
By (iv) and (v), ${|K|_p=|K:P_n\cap K|}_p={|P_nK:P_n|}_p\leq p^n$. If $K\neq 1$, then $N\leq K$, and so $G/K$ is $p$-nilpotent. It follows from Lemma 2.12 that $G$ is $p$-nilpotent, a contradiction. We may, therefore, assume that $K=1$. Then $P_n$ is S-quasinormal in $G$. If $P_n\neq 1$, then $N\leq {P_n}^G\leq P_1$ as above, a contradiction. Hence $P_n=1$, and thus $|G|_p=p^n$, the final contradiction completes the proof of (7).\par
(8) \textit{$H=G$.}\par
If not, then $H$ is $p$-nilpotent by (4). Let $H_{p'}$ be the normal $p$-complement of $H$. Then $H_{p'}\unlhd G$, which induces that $H=O_p(H)=1$ by (6) and (7). So $G$ is $p$-nilpotent, a contradiction.\par
(9) \textit{$O_p(G)=1$.}\par
It follows directly from (7) and (8).\par
(10) \textit{$N$ is not $p$-soluble.}\par
If $N$ is $p$-soluble, then $O_{p'}(N)\neq 1$ or $O_p(N)\neq 1$, which contradicts (6) or (9).\par
(11) \textit{Let $P$ be a Sylow $p$-subgroup of $G$, then $G=PN$ and $P\cap N\nleq \Phi(P)$. In addition, $N=O^p(G)$.}\par
If $PN<G$, then $PN$ is $p$-nilpotent by (4) and so is $N$, contrary to (10). Hence $G=PN$. If $P\cap N\leq \Phi(P)$, then $N$ is $p$-nilpotent by Lemma 2.17, which also contradicts (10). Clearly $1\neq O^p(G)\leq N$, and so $N=O^p(G)$.\par
(12) \textit{Final contradiction.}\par
By (11), there exists a maximal subgroup $P_1$ of $P$ such that $P=(P\cap N)P_1$. It is obvious that $P\cap G^{\mathfrak{N}_\mathfrak{p}}=P\cap N\nleq P_1$. Let $P_n$ be an $n$-maximal subgroup of $P$ contained in $P_1$. By the hypothesis, $P_n$ either satisfies $(\triangle)$ or has a $p$-nilpotent supplement $T$ in $G$.\par
(i) \textit{Case \textup{1:} $P_n$ satisfies $(\triangle)$ in $G$.}\par
In this case, there exists a normal subgroup $K$ of $G$ such that $P_nK$ is S-quasinormal in $G$ and $P_n\cap K\leq (P_n)_{seG}$ or $P_n\cap K\leq (P_n)_{\tau G}$. If $K=1$, then $P_n$ is S-quasinormal in $G$. It induces from Lemma 2.1(1) and (9) that $P_n\leq O_p(G)=1$, and so $|G|_p\leq p^n$, which contradicts (1). Therefore we may suppose that $K\neq 1$. Then $N\leq K$, and so $P_n\cap N\leq (P_n)_{seG}$ or $P_n\cap N\leq (P_n)_{\tau G}$.\par
First assume that $P_n\cap N\leq (P_n)_{seG}$. With the similar argument as above, let $D_1,D_2,\cdots D_s$ be all subgroups of $P_n$ which are S-quasinormally embedded in $G$. Then there exist S-quasinormal subgroups $U_1,U_2,\cdots U_s$ of $G$ with $D_i\in Syl_p(U_i)$ $(1\!\leq\! i \!\leq \! s)$. If $(U_i)_G\neq 1$, then $N\leq (U_i)_G\leq U_i$. Thus $D_i\cap N\in Syl_p(N)$. Since $P\cap N\in Syl_p(N)$, $D_i\cap N=P\cap N$. It follows that $P\cap N\leq D_i\leq P_n\leq P_1$, a contradiction. Therefore $(U_i)_G=1$, and so $D_i$ is S-quasinormal in $G$ by Lemma 2.3. As $(P_n)_{seG}=\langle D_1,D_2,\cdots D_s \rangle$, $(P_n)_{seG}\leq O_p(G)=1$, which induces that $P_n\cap N=1$. Hence $|N|_p=|P\cap N|=|P\cap N:P_n\cap N|=|P_n(P\cap N):P_n|\leq |P:P_n|=p^n$. It follows from Lemma 2.12 that $G$ is $p$-nilpotent, a contradiction.\par
Now consider that $P_n\cap N\leq (P_n)_{\tau G}$. Let $q\neq p$ be an arbitrary prime divisor of $|G|$ and $Q\in Syl_q(G)$. Clearly $Q\leq O^p(G)=N$ and $(p,|Q^G|)\neq 1$. Then by Lemma 2.7(1), $(P_n)_{\tau G}Q=Q(P_n)_{\tau G}$. It follows that $(P_n\cap N)Q=((P_n)_{\tau G}\cap N)Q=(P_n)_{\tau G}Q\cap N=Q(P_n)_{\tau G}\cap N=Q((P_n)_{\tau G}\cap N)=Q(P_n\cap N)$. Therefore, $(P_n\cap N)Q^n=Q^n(P_n\cap N)$ for all $n\in N$. If $N=(P_n\cap N)Q$, then $N$ is soluble by Burnside's $p^aq^b$-Theorem, which contradicts (10). So there exists a non-trivial normal subgroup $X$ of $N$ such that either $P_n\leq X$ or $Q\leq X$ by Lemma 2.13. Since $N$ is characteristically simple group of $G$, $N\cong A_1\times A_2\times \cdots \times A_t$, where $A_i\cong A$ $(1\!\leq\! i \!\leq \! t)$ is a simple group. Obviously $A$ is non-abelian owing to (10). Without loss of generality, we may assume that $X\cong A_1\times A_2\times \cdots \times A_k$ $(k<t)$ by Lemma 2.18. Since $Q\in Syl_q(N)$, $|Q|=|N|_q={(|A|_q)}^t>{(|A|_q)}^k=|X|_q$. This implies that $P_n\cap N\leq X$. Note that $|P\cap N:P_n\cap N|=|P_n(P\cap N):P_n|\leq |P:P_n|=p^n$, then we have $|N|_p/|X|_p={(|A|_p)}^{t-k}\leq p^n$. As $t-k\geq 1$, $|A|_p\leq p^n$. Hence $A$ is $p$-nilpotent by Lemma 2.12. Consequently, $N$ is $p$-soluble, a contradiction.\par
(ii) \textit{Case \textup{2:} $P_n$ has a $p$-nilpotent supplement $T$ in $G$.}\par
Let $T_{p'}$ be the normal $p$-complement of $T$. Then $G=P_nT=P_nN_G(T_{p'})$ and $N_G(T_{p'})<G$. Since $P=P_n(P\cap N_G(T_{p'}))$ and $P\cap N_G(T_{p'})<P$, there exists a maximal subgroup $P_2$ of $P$ such that $P\cap N_G(T_{p'})\leq P_2$. As $N=O^p(G)$ by (11), every Hall $p'$-subgroup of $G$ is contained in $N$. By Lemma 2.14 and Frattni argument, $G=NN_G(T_{p'})=PN_G(T_{p'})$. Then $|P|=|G|_p\leq |N|_p|N_G(T_{p'})|_p=|P\cap N||P\cap N_G(T_{p'})|$. Thus $P=(P\cap N)(P\cap N_G(T_{p'}))\leq (P\cap N)P_2$. This implies that $P=(P\cap N)P_2$ and $P\cap G^{\mathfrak{N}_\mathfrak{p}}=P\cap N\nleq P_2$.\par
Let $P_{n2}$ be an $n$-maximal subgroup of $P$ contained in $P_2$. Then by the hypothesis, $P_{n2}$ either has a $p$-nilpotent supplement $E$ or satisfies $(\triangle)$ in $G$. If $P_{n2}$ has a $p$-nilpotent supplement $E$ in $G$, with the similar discussion as the proof of (2)(ii), we can obtain a contradiction. Hence $P_{n2}$ satisfies $(\triangle)$ in $G$. Then the same argument as (i) yields a contradiction again. The theorem is proved.\par
\medskip
\noindent\textbf{Theorem 3.2.} \textit{Let G be a group and p a prime divisor of $|G|$ with $(|G|,(p-1)(p^2-1)\cdots (p^n-1))=1$ for some integer $n\geq 1$. Then G is p-nilpotent if and only if there exists a normal subgroup H of G such that G/H is p-nilpotent and for any Sylow p-subgroup P of H, every subgroup L of $P\cap G^{\mathfrak{N}_\mathfrak{p}}$ of order $p^n$ or order \textup{4} \textup{(}when $p=2$, $n=1$, P is non-abelian, and $L$ is cyclic\textup{)} not contained in $Z_{\infty}(G)$ \textup{(}if exists\textup{)} either has a p-nilpotent supplement in G or satisfies $(\triangle)$ in $G$.}\par
\medskip
\noindent\textit{Proof.} We need only to prove the sufficiency. Suppose that the result is false and let $(G,H)$ be a counterexample for which $|G|+|H|$ is minimal. Then:\par
(1) \textit{${|H|}_p\geq p^{n+1}$.}\par
It follows from Lemma 2.12.\par
(2) \textit{If $1\neq N\unlhd G$ with $(|N|,p)=1$, then $G/N$ is $p$-nilpotent.}\par
If $|(HN/N)\cap (G/N)^{\mathfrak{N}_\mathfrak{p}}|_p=|G^{\mathfrak{N}_\mathfrak{p}}N/N|_p\leq p^n$, then $G/N$ is $p$-nilpotent by Lemma 2.12. So we may assume that $|G^{\mathfrak{N}_\mathfrak{p}}N/N|_p\geq p^{n+1}$. Let $L/N$ be a subgroup of $PN/N\cap G^{\mathfrak{N}_\mathfrak{p}}N/N$ of order $p^n$ or $4$ (when $p=2$, $n=1$, $PN/N$ is non-abelian, and $L/N$ is cyclic) not contained in $Z_{\infty}(G/N)$, where $PN/N\in Syl_p(HN/N)$ and $P\in Syl_p(H)$. Since $L=(P\cap L)N$ and $(|N|,p)=1$, $|L/N|=|(P\cap L)N/N|=|P\cap L|=p^n$ or 4. As $P\cap L\leq G^{\mathfrak{N}_\mathfrak{p}}N$ and $(|P\cap L|,|G^{\mathfrak{N}_\mathfrak{p}}N:G^{\mathfrak{N}_\mathfrak{p}}|)=1$, $P\cap L\leq P\cap G^{\mathfrak{N}_\mathfrak{p}}$. By Lemma 2.19(1), $P\cap L\nsubseteq Z_{\infty}(G)$. Suppose that $|P\cap L|=4$, then $P$ is non-abelian, and $P\cap L$ is cyclic owing to $G$-isomorphism $L/N\cong P\cap L$. Hence by the hypothesis, it is easy to see that $L/N$ either has a $p$-nilpotent supplement or satisfies $(\triangle)$ in $G/N$ with the similar discussion as the step (3) in the proof of Theorem 3.1. This means that $(G/N,HN/N)$ satisfies the hypothesis. The choice of $(G,H)$ implies that $G/N$ is $p$-nilpotent.\par
(3) \textit{Every non-trivial subgroup $F$ of $G$ is $p$-nilpotent.}\par
Considering $(F,H\cap F)$, we may assume that $|(H\cap F)\cap F^{\mathfrak{N}_\mathfrak{p}}|_p=|F^{\mathfrak{N}_\mathfrak{p}}|_p\geq p^{n+1}$. If not, then $F$ is $p$-nilpotent by Lemma 2.12. Let $L$ be a subgroup of $(P\cap F)\cap F^{\mathfrak{N}_\mathfrak{p}}=P\cap F^{\mathfrak{N}_\mathfrak{p}}$ of order $p^n$ or 4 (when $p=2$, $n=1$, $P\cap F$ is non-abelian, and $L$ is cyclic) not contained in $Z_{\infty}(F)$, where $P\cap F\in Syl_p(H\cap F)$ and $P\in Syl_p(H)$. By Lemma 2.19(2), $L\nsubseteq Z_{\infty}(G)$. In view of that $P\cap F^{\mathfrak{N}_\mathfrak{p}}\leq P\cap (F \cap G^{\mathfrak{N}_\mathfrak{p}}) \leq P\cap G^{\mathfrak{N}_\mathfrak{p}}$, $L$ either has a $p$-nilpotent supplement or satisfies $(\triangle)$ in $F$ by Lemma 2.9(1) and Lemma 2.10(1). Hence $(F,H\cap F)$ satisfies the hypothesis. Then $F$ is $p$-nilpotent by the choice of $(G,H)$.\par
(4) \textit{$O_{p'}(G)=1$.}\par
If not, by (2), $G/O_{p'}(G)$ is $p$-nilpotent and so is $G$, a contradiction.\par
(5) \textit{$G$ is a minimal nonnilpotent group.}\par
By (3), $G$ is a minimal non-$p$-nilpotent group. Then $G$ is a minimal nonnilpotent group by It\v{o}' s theorem (see \cite[Chapter IV, Theorem 5.4]{Hup}). Therefore by \cite[Theorem 3.4.11]{Guo4} and Lemma 2.11, $G=P\rtimes Q$, where $P$ is a normal Sylow $p$-subgroup of $G$, and $Q$ is a cyclic Sylow $q$-subgroup of $G$ ($p\neq q$); and:\par (i) $P/\Phi(P)$ is a chief factor of $G$.\par (ii) $P=G^{\mathfrak{N}}=G^{\mathfrak{N}_\mathfrak{p}}$.\par (iii) $exp(P)=p$ or 4 (when $p=2$, $P$ is non-abelian).\par (iv) $\Phi(G)=Z_{\infty}(G)$.\par (v) $F(G)=F_p(G)=P\Phi(G)$.\par
(6) \textit{$P\leq H$.}\par
It follows from the fact that $P=G^{\mathfrak{N}_\mathfrak{p}}$ and $G/H$ is $p$-nilpotent.\par
(7) \textit{$F(G)=P$ and $\Phi(G)=\Phi(P)$.}\par
Since $F_p(G)/O_{p'}(G)=O_p(G/O_{p'}(G))$, by (4) and (5), $F(G)=F_p(G)=P\Phi(G)=O_p(G)=P$. Thus $\Phi(P)\leq \Phi(G)\leq P$. As $P/\Phi(P)$ is a chief factor of $G$, if $\Phi(G)=P$, then $G=Q$, a contradiction. So we have $\Phi(G)=\Phi(P)$.\par
(8) \textit{$P$ has a proper subgroup $L$ of order $p^n$ or $4$ such that $L\nleq \Phi(P)$ and $L$ either has a $p$-nilpotent supplement or satisfies $(\triangle)$ in $G$.}\par
Take an element $x\in P\backslash \Phi(P)$, and let $E=\langle x\rangle$. Then $|E|=p$ or 4 (when $p=2$, $P$ is non-abelian). Clearly, $E\nsubseteq Z_{\infty}(G)=\Phi(G)=\Phi(P)$. It induces from Lemma 2.20 that there exists a subgroup $L$ of $P$ order $p^n$ or 4 (when $p=2$, $n=1$, we may take $L=E$) such that $E\leq L$ and $L\nleq \Phi(P)$. By the hypothesis, $L$ either has a $p$-nilpotent supplement $T$ or satisfies $(\triangle)$ in $G$. Moreover, if $L=P$, then since $|P|\geq p^{n+1}$ by (1), $|P|=4$. This implies that $P$ is abelian, a contradiction.\par
(9) \textit{Final contradiction.}\par
(i) \textit{Case \textup{1:} $L$ has a $p$-nilpotent supplement $T$ in $G$.}\par
Since $G=LT=PT$, $P\cap T\unlhd T$ and $T\Phi(P)/\Phi(P)\leq N_{G/\Phi(P)}((P\cap T)\Phi(P)/\Phi(P))$. Obviously, $(P\cap T)\Phi(P)/\Phi(P)\unlhd P/\Phi(P)$ for $P/\Phi(P)$ is abelian. Hence $(P\cap T)\Phi(P)/\Phi(P)\unlhd G/\Phi(P)$. As $P/\Phi(P)$ is a chief factor of $G$, $(P\cap T)\Phi(P)=\Phi(P)$ or $P$. In the former case, $P\cap T\leq \Phi(P)$. This implies that $P=L(P\cap T)=L$, which contradicts (8). In the latter case, $P\cap T=P$. Then $P\leq T$, and so $G=T$ is $p$-nilpotent, also a contradiction.\par
(ii) \textit{Case \textup{2:} $L$ satisfies $(\triangle)$ in $G$.}\par
By the hypothesis, there exists a normal subgroup $K$ of $G$ such that $LK$ is S-quasinormal in $G$ and $L\cap K\leq L_{seG}$ or $L\cap K\leq L_{\tau G}$. Since $L\leq P=O_p(G)$, by Lemma 2.8 and Lemma 2.2, $L_{seG}=L_{\tau G}=L_{sG}$, where $L_{sG}$ is the largest S-quasinormal subgroup of $G$ contained in $L$.\par
By Lemma 2.1(5), $O^p(G)\leq N_G(L_{sG})$. Hence $O^p(G)\Phi(P)/\Phi(P)\leq N_G(L_{sG})\Phi(P)/\Phi(P)\leq N_{G/\Phi(P)}(L_{sG}\Phi(P)/\Phi(P))$. As $P/\Phi(P)$ is abelian, $P/\Phi(P)$ $\leq N_{G/\Phi(P)}(L_{sG}\Phi(P)/\Phi(P))$. It follows that $L_{sG}\Phi(P)/\Phi(P)\unlhd G/\Phi(P)$, and so $L_{sG}\Phi(P)\unlhd G$. Therefore $L_{sG}\Phi(P)=P$ or $\Phi(P)$.\par
If $L_{sG}\Phi(P)=P$, then $L_{sG}=P=L$, which contradicts (8). Now assume that $L_{sG}\Phi(P)=\Phi(P)$, then $L_{sG}\leq \Phi(P)$. If $K=G$, then $L=L_{sG}\leq \Phi(P)$, a contradiction. Thereupon $K<G$. Since $K$ is nilpotent, $K\leq F(G)=P$. This shows that $LK$ is an S-quasinormal subgroup of $G$ contained in $P$. With the similar discussion as above, $LK\Phi(P)\unlhd G$, and so $LK\Phi(P)=P$ or $\Phi(P)$. If $LK\Phi(P)=\Phi(P)$, then $L\leq LK\leq \Phi(P)$, a contradiction. Hence $LK\Phi(P)=P$, then $LK=P$. As $K\Phi(P)\unlhd G$, $K\Phi(P)=P$ or $\Phi(P)$. In the former case, $K=P$, and thereby $L=L_{sG}\leq \Phi(P)$, a contradiction. In the latter case, $K\leq \Phi(P)$, and thus $L=P$, a contradiction again. The theorem is proved.\par
\medskip
\noindent\textbf{Theorem 3.3.} \textit{Let G be an $A_4$-free group and p a prime divisor of $|G|$ with $(|G|,p-1)=1$. Then G is p-nilpotent if and only if there exists a normal subgroup H of G such that G/H is p-nilpotent and for any Sylow p-subgroup P of H, every 2-maximal subgroup of P not containing $P\cap G^{\mathfrak{N}_\mathfrak{p}}$ or every subgroup of order $p^2$ of $P\cap G^{\mathfrak{N}_\mathfrak{p}}$ not contained in $Z_{\infty}(G)$ \textup{(}if exists\textup{)} either has a p-nilpotent supplement in G or satisfies $(\triangle)$ in $G$.}\par
\medskip
\noindent\textit{Proof.} In view of Lemma 2.21, the proof is analogous to Theorem 3.1 and 3.2.\par
\medskip
\noindent\textbf{Theorem 3.4.} \textit{Let $\mathfrak{F}$ be a saturated formation containing $\mathfrak{U}$ and $G$ a group. Then $G\in \mathfrak{F}$ if and only if there exists a normal subgroup H of G such that $G/H\in \mathfrak{F}$ and for any non-cyclic Sylow subgroup P of H, every maximal subgroup of P either has a supersolvable supplement in G or satisfies $(\triangle)$ in G.}\par
\medskip
\noindent\textit{Proof.} We need only to prove the sufficiency. Assume that the result is false and take $(G,H)$ a counterexample for which $|G|+|H|$ is minimal. Then:\par
(1) \textit{If $r$ is the maximal prime divisor of $|H|$, then the Sylow $r$-subgroup $R$ of $H$ is normal in $G$.}\par
Let $p$ be the minimal prime divisor of $|H|$ and $P$ a Sylow $p$-subgroup of $H$. If $P$ is cyclic, then $H$ is $p$-nilpotent by Lemma 2.22. Note that a supersolvable group is $p$-nilpotent when $p$ is the smallest prime divisor of its order. If $P$ is non-cyclic, then in view of that $(H,H)$ satisfies the hypothesis, we also have $H$ is $p$-nilpotent by Theorem 3.1. Therefore $H$ has the normal $p$-complement $H_{p'}$ such that $H_{p'}\,char\,H\unlhd G$. Clearly, $(H_{p'},H_{p'})$ also satisfies the hypothesis. Let $q$ be the minimal prime divisor of $|H_{p'}|$, then it follows that $H_{p'}$ is $q$-nilpotent as above. The rest may be deduced by analogy. Thus, if $r$ is the maximal prime divisor of $|H|$ and $R$ is a Sylow $r$-subgroup of $H$, we obtain that $R\unlhd G$.\par
(2) \textit{There exists a unique minimal normal subgroup $N$ of $G$ contained in $R$ such that $G/N\in \mathfrak{F}$ and $N\nleq \Phi(G)$.}\par
Indeed, it is easy to see that $(G/N,H/N)$ satisfies the hypothesis for any minimal normal subgroup $N$ of $G$ contained in $R$. The choice of $(G,H)$ implies that $G/N\in \mathfrak{F}$. Since $\mathfrak{F}$ is a saturated formation, $N$ is the unique minimal normal subgroup of $G$ contained in $R$ and $N\nleq \Phi(G)$.\par
(3) \textit{$N=R\leq O_r(G)$.}\par
Since $N\nleq \Phi(G)$, there exists a maximal subgroup $M$ of $G$ such that $G=NM=RM$. By Lemma 2.16, $R\cap M$ is normal in $G$. Then $R\cap M=1$ for $N\nleq R\cap M$. Hence $R=N(R\cap M)=N\leq O_r(G)$.\par
(4) \textit{Final contradiction.}\par
If $R$ is cyclic, then $G\in \mathfrak{F}$ by Lemma 2.23, a contradiction. Thus $R$ is non-cyclic. Let $R'$ be a Sylow $r$-subgroup of $G$, then $R'=R(M\cap R')$. We may /choose a maximal subgroup $R_1'$ of $R'$ such that $M\cap R'\leq R_1'$. Clearly, $R_1=R\cap R_1'$ is a maximal subgroup of $R$ with $R_1\unlhd R'$. By the hypothesis, $R_1$ either has a supersolvable supplement $T$ or satisfies $(\triangle)$ in $G$.\par
First suppose that $G=R_1T=RT$, where $T$ is a supersolvable subgroup of $G$. Then since $R=N$ is abelian, $R\cap T\unlhd G$. It follows that $R\cap T=1$ or $R$ . In the former case, $R=R_1(R\cap T)=R_1$, a contradiction. In the latter case, $R\leq T$, and thereby $G=T\in \mathfrak{F}$, also a contradiction.\par
Now assume that $R_1$ satisfies $(\triangle)$ in $G$. Then there exists a normal subgroup $K$ of $G$ such that $R_1K$ is S-quasinormal in $G$ and $R_1\cap K\leq (R_1)_{seG}$ or $R_1\cap K\leq (R_1)_{\tau G}$. Since $R_1\leq O_r(G)$, $(R_1)_{seG}$ and $(R_1)_{\tau G}$ are both S-quasinormal in $G$ by Lemma 2.8 and Lemma 2.2. Hence $(R_1)_{seG}=(R_1)_{\tau G}=(R_1)_{sG}$, where $(R_1)_{sG}$ is the largest S-quasinormal subgroup of $G$ contained in $R_1$. Obviously $R\cap K=1$ or $R$. In the former case, $R_1=R\cap R_1K$ is S-quasinormal in $G$ by Lemma 2.2. In the latter case, $R\leq K$, and so $R_1=R_1\cap K=(R_1)_{sG}$ is also S-quasinormal in $G$. However, since $R_1\unlhd R'$ and $O^r(G)\leq N_G(R_1)$, $R_1\unlhd G$. Therefore $R_1=1$. This implies that $R$ is cyclic, and consequently $G\in \mathfrak{F}$ by Lemma 2.23. The final contradiction completes the proof.\par
\medskip
\noindent\textbf{Theorem 3.5.} \textit{Let $\mathfrak{F}$ be a saturated formation containing $\mathfrak{U}$ and $G$ a group. Then $G\in \mathfrak{F}$ if and only if there exists a normal subgroup H of G such that $G/H\in \mathfrak{F}$ and for any non-cyclic Sylow subgroup P of H, every cyclic subgroup L of P of prime order or order \textup{4} \textup{(}when $p=2$, P is non-abelian\textup{)} not contained in $Z_{\infty}(G)$ either has a supersolvable supplement in G or satisfies $(\triangle)$ in G.}\par
\medskip
\noindent\textit{Proof.} We need only to prove the sufficiency. Assume that the result is false and take $(G,H)$ a counterexample for which $|G|+|H|$ is minimal. Then:\par
(1) \textit{The Sylow $q$-subgroup $Q$ of $H$ is normal in $G$, where $q$ is the maximal prime divisor of $|H|$.}\par
In view of Lemma 2.22 and Theorem 3.2, we can obtain the assertion as the step (1) in the proof of Theorem 3.4.\par
(2) \textit{$G^{\mathfrak{F}}=Q=H$.}\par
By Lemma 2.9(3) and Lemma 2.10(3), it is easy to see that $(G/Q,H/Q)$ satisfies the hypothesis. The choice of $(G,H)$ implies that $G/Q\in \mathfrak{F}$. Hence $G^{\mathfrak{F}}\leq Q\leq H$. Since $1\neq G^{\mathfrak{F}}$ is a $q$-subgroup of $H$, $(G,G^{\mathfrak{F}})$ also satisfies the hypothesis. This induces that $G^{\mathfrak{F}}=Q=H$.\par
(3) \textit{Let $M$ be any maximal subgroup of $G$ not containing $G^{\mathfrak{F}}=Q$, then $M\in \mathfrak{F}$.}\par
In fact, $G=MQ$ and $M/M\cap Q\cong MQ/Q=G/Q\in \mathfrak{F}$. By Lemma 2.9(1) and Lemma 2.10(1), $(M,M\cap Q)$ satisfies the hypothesis. Thus $M\in \mathfrak{F}$ by the choice of $(G,H)$.\par
(4) \textit{$Q/\Phi(Q)$ is a chief factor of $G$ and $exp(Q)=q$ or $4$ \textup{(}when $q=2$, $Q$ is non-abelian\textup{)}.}\par
It follows directly from (3) and \cite[Theorem 3.4.2]{Guo4}.\par
(5) \textit{$Q$ has a subgroup $L$ of order $q$ or $4$ such that $L\nleq \Phi(Q)$ and $L$ either has a supersolvable supplement or satisfies $(\triangle)$ in $G$.}\par
Take an element $x\in Q\backslash \Phi(Q)$, and let $L=\langle x\rangle$. Then $|L|=q$ or 4 (when $q=2$, $Q$ is non-abelian). Assume that $L\subseteq Z_{\infty}(G)$. By (4), $(Q\cap Z_{\infty}(G))\Phi(Q)=\Phi(Q)$ or $Q$. In the former case, $L\leq Q\cap Z_{\infty}(G)\leq \Phi(Q)$, a contradiction. In the latter case, $Q\leq Z_{\infty}(G)\leq Z_{\mathfrak{U}}(G)$, where $Z_{\mathfrak{U}}(G)$ is the product of all such normal subgroups $E$ of $G$ whose $G$-chief factors have prime order. Hence $|Q/\Phi(Q)|=q$. This induces that $Q$ is cyclic, and thus $G\in \mathfrak{F}$ by Lemma 2.23, a contradiction. Therefore, $L\nsubseteq Z_{\infty}(G)$. Then by the hypothesis, $L$ either has a supersolvable supplement $T$ or satisfies $(\triangle)$ in $G$.\par
(6) \textit{$L\neq Q$.}\par
If not, then $Q=G^{\mathfrak{F}}$ is cyclic. By Lemma 2.23, $G\in \mathfrak{F}$, a contradiction.\par
(7) \textit{Final contradiction.}\par
First suppose that $L$ satisfies $(\triangle)$ in $G$. With the similar argument as the step (9) in the proof of Theorem 3.2, we see that it is impossible. Now assume that $G=LT$ and $T$ is supersolvable. Let $r$ be the maximal prime divisor of $|T|$ and $R$ a Sylow $r$-subgroup of $T$. As $L\neq Q$, $r\geq q$. Since $R\unlhd T$, $G=LN_G(R)$ and $|G:N_G(R)|=1$ or $q$ or 4.\par
(i) \textit{Case \textup{1:} $|G:N_G(R)|=1$.}\par
In this case, $R\unlhd G$. If $r>q$, then $R$ is a Sylow $r$-subgroup of $G$. As $(G/R,QR/R)$ satisfies the hypothesis, $G/R\in \mathfrak{F}$ by the choice of $(G,H)$. This implies that $G^{\mathfrak{F}}=Q\leq R$, a contradiction. Now assume that $r=q$. Then $(R\cap Q)\Phi(Q)=\Phi(Q)$ or $Q$ for $Q/\Phi(Q)$ is a chief factor of $G$. In the former case, $R\cap Q\leq \Phi(Q)$. Since $Q\cap T\leq R$, $Q=L(Q\cap T)\leq LR$. Thus $Q=LR\cap Q=L(R\cap Q)$, and so $Q=L$, which contradicts (6). In the latter case, $L\leq Q\leq R \leq T$. Then $G=T\in {\mathfrak{F}}$, also a contradiction.\par
(ii) \textit{Case \textup{2:} $|G:N_G(R)|=2$.}\par
In this case, $N_G(R)\unlhd G$. This implies that $Q=G^{\mathfrak{F}}\leq N_G(R)$, and thereby $G=N_G(R)$, a contradiction.\par
(iii) \textit{Case \textup{3:} $|G:N_G(R)|=q$, where $q$ is an odd prime.}\par
In this case, $L\cap N_G(R)=1$, and so $N_G(R)=(L\cap N_G(R))T=T$. If $r>q$, then $R$ is a Sylow $r$-subgroup of $G$. As $|G:N_G(R)|=q$, the number of the Sylow $r$-subgroups is $q$. But by Sylow's Theorem, $q\equiv 1$ (mod $r$), which is impossible. Now consider that $r=q$. Then $T\leq N_G(Q\cap R)$. Since $T$ is a maximal subgroup of $G$, $N_G(Q\cap R)=G$ or $T$. In the former case, $Q\cap R\unlhd G$, with the similar discussion as (i), we can obtain a contradiction. In the latter case, clearly $Q\cap R<Q$. This shows that $Q\cap R<N_Q(Q\cap R)=Q\cap N_G(Q\cap R)=Q\cap T$. But $Q\cap T\leq R$, also a contradiction.\par
(iv) \textit{Case \textup{4:} $|G:N_G(R)|=4$.}\par
In this case, $q=2$ and $T=N_G(R)$. If $r>q=2$, then since $R$ is a Sylow $r$-subgroup of $G$, $4\equiv 1$ (mod $r$). Hence $r=3$, and thus $\pi(G)=\{2,3\}$. As $|L|=4$, $Q$ is non-cyclic by (4). Note that $G/Q=QT/Q\cong T/Q\cap T$ is supersolvable and so 2-nilpotent. We obtain that $(G,Q)$ satisfies the hypothesis of Theorem 3.2 (take $p=2$, $n=1$). Consequently, $G$ is 2-nilpotent. This implies that $G_{2'}=R\unlhd G$, and thereby $T=N_G(R)=G$, a contradiction. Now suppose that $r=q=2$. Then $T\leq N_G(Q\cap R)$. This induces that $|G:N_G(Q\cap R)|=1$ or 2 or 4. If $|G:N_G(Q\cap R)|=1$ or 4, then $N_G(Q\cap R)=G$ or $T$, a contradiction can be obtained as (iii). Therefore $|G:N_G(Q\cap R)|=2$. Then $N_G(Q\cap R)\unlhd G$, and thus $Q=G^{\mathfrak{F}}\leq N_G(Q\cap R)$. Hence $N_G(Q\cap R)=G$, the final contradiction completes the proof.\par
\medskip
\noindent\textbf{Theorem 3.6.} \textit{Let $\mathfrak{F}$ be a saturated formation containing $\mathfrak{U}$. Then $G\in \mathfrak{F}$ if and only if there exists a normal subgroup H of G such that $G/H\in \mathfrak{F}$ and for any non-cyclic Sylow subgroup P of $F^*(H)$, at least one of the following holds:}\par
\textit{\textup{(i)} Every maximal subgroup of P not having a supersolvable supplement satisfies $(\triangle)$ in G.}\par
\textit{\textup{(ii)} Every cyclic subgroup $L$ of P of prime order or order \textup{4} $\bigm{(}$when $p=2$, P is non-abelian, and $L\nsubseteq  Z_{\infty}(G)$$\bigm{)}$ not having a supersolvable supplement satisfies $(\triangle)$ in G.}\par
\medskip
\noindent\textit{Proof.} We need only to prove the sufficiency.\par
(1) \textit{\textit{Case \textup{1:}} $\mathfrak{F}=\mathfrak{U}$.}\par
Obviously, $\bigm{(}$$F^*(H)$,$\:$$F^*(H)$$\bigm{)}$ satisfies the hypothesis of Theorem 3.4 or 3.5. Hence $F^*(H)$ is supersolvable and so soluble. It induces from \cite[Chapter X, Corollary 13.7(d)]{Hup1} that $F^*(H)=F(H)\leq F(G)$. In this case, weakly $\tau$-embedding is equivalent to weakly S-embedding by Lemma 2.8. It follows that $(G,H)$ satisfies the hypothesis of Lemma 2.24. Thereupon $G\in \mathfrak{F}$.\par
(2) \textit{\textit{Case \textup{2:}} $\mathfrak{F}\neq \mathfrak{U}$.}\par
Obviously, $(H,H)$ satisfies the hypothesis of (1). Hence $H$ is supersolvable and so soluble. This implies that $F^*(H)=F(H)\leq F(G)$. Therefore $G\in \mathfrak{F}$ by Lemma 2.24 again.\par

\section{Applications}
In Section 1, we have seen that every S-quasinormal, SS-quasinormal,
S-semipermutable, $\tau$-quasinormal, and S-embedded subgroup of $G$
are all weakly $\tau$-embedded in $G$. On the other hand, it is
clear that every S-quasinormally embedded subgroup of $G$ is weakly
S-embedded in $G$.\par Except the concepts mentioned in Section 1,
there are several other notions were introduced. Recall that: a
subgroup $H$ of $G$ is said to be c-normal \cite{Wan} in $G$ if
there exists a normal subgroup $K$ of $G$ such that $G=HK$ and
$H\cap K\leq H_G$, where $H_G$ is the largest normal subgroup of $G$
contained in $H$; a subgroup $H$ of $G$ is said to be
$\mathrm{c^*}$-normal \cite{Wei} in $G$ if there exists a normal
subgroup $K$ of $G$ such that $G=HK$ and $H\cap K$ is
S-quasinormally embedded in $G$; a subgroup $H$ of $G$ is said to be
n-embedded \cite{Guo2} in $G$ if there exists a normal subgroup $K$
of $G$ such that $HK=H^G$ and $H\cap K\leq H_{sG}$. It is easy to
see that all the above subgroups of $G$ are also either weakly
S-embedded or weakly $\tau$-embedded in $G$.\par Therefore, a large
number of results are deduced immediately from our theorems. Here we
only list some of recent ones.\par \noindent\textbf{Corollary 4.1}
\cite[Theorem 3.3]{Wan2}. \textit{Let p be the smallest prime number
dividing the order of G and P a Sylow p-subgroup of G. If every
maximal subgroup of P is S-semipermutable in G, then G is
p-nilpotent.}\par \noindent\textbf{Corollary 4.2} \cite[Theorem
3.1]{Wei}. \textit{Let H be a normal subgroup of G such that G/H is
p-nilpotent and P a Sylow p-subgroup of H, where p is a prime
divisor of $|G|$ with $(|G|,p-1)=1$. If all maximal subgroups of P
are $c^*$-normal in G, then G is p-nilpotent. In particular, G is
p-supersolvable.}\par \noindent\textbf{Corollary 4.3} \cite[Theorem
2.3]{Guo3}. \textit{Let P be a Sylow p-subgroup of G, where p is a
prime divisor of $|G|$ with $(|G|,p-1)=1$. If every maximal subgroup
of P not having a p-nilpotent supplement in G is S-embedded in G,
then G is p-nilpotent.}\par \noindent\textbf{Corollary 4.4}
\cite[Theorem 1.7]{Li}. \textit{Let p be the smallest prime dividing
the order of G and P a Sylow p-subgroup of G. If every
\textup{2}-maximal subgroup of P is SS-quasinormal in G and G is
$A_4$-free, then G is p-nilpotent.}\par \noindent\textbf{Corollary
4.5} \cite[Theorem 3.5]{Wan2}. \textit{Let p be the smallest prime
number dividing the order of G and P is a Sylow p-subgroup of G. If
every \textup{2}-maximal subgroup of P is S-semipermutable in G and
G is $A_4$-free, then G is p-nilpotent.}\par
\noindent\textbf{Corollary 4.6} \cite[Theorem 2]{Zha}. \textit{Let
$\mathfrak{F}$ be a saturated formation containing $\mathfrak{U}$.
Then $G\in \mathfrak{F}$ if and only if there exists a normal
subgroup H of G such that $G/H\in \mathfrak{F}$ and all maximal
subgroups of any Sylow subgroup of H are S-semipermutable in G.}\par
\noindent\textbf{Corollary 4.7} \cite[Theorem 4.1]{Wei}. \textit{Let
$\mathfrak{F}$ be a saturated formation containing $\mathfrak{U}$.
Suppose that $G$ is a group with a normal subgroup H such that
$G/H\in \mathfrak{F}$. If all maximal subgroups of any Sylow
subgroup of H are $c^*$-normal in G, then $G\in \mathfrak{F}$.}\par
\noindent\textbf{Corollary 4.8} \cite[Theorem 3.4]{Li2}. \textit{If
all cyclic subgroups of prime order or order \textup{4} of G are
SS-quasinormal in G, then G is supersolvable.}\par
\noindent\textbf{Corollary 4.9} \cite[Theorem 3]{Zha}. \textit{Let
$\mathfrak{F}$ be a saturated formation containing $\mathfrak{U}$.
Then $G\in \mathfrak{F}$ if and only if there exists a normal
subgroup H of G such that $G/H\in \mathfrak{F}$ and all cyclic
subgroups of prime order or order \textup{4} of H are
S-semipermutable in G.}\par \noindent\textbf{Corollary 4.10}
\cite[Theorem D]{Guo2}. \textit{Let $\mathfrak{F}$ be a saturated
formation containing $\mathfrak{U}$ and $G$ a group with a normal
subgroup E such that $G/E\in \mathfrak{F}$. If for every non-cyclic
Sylow subgroup P of E, every maximal subgroup of P or every cyclic
subgroup H of P of prime order or order \textup{4} $\bigm{(}$if P is
non-abelian \textup{2}-group and $H\nsubseteq
Z_{\infty}(G)$$\bigm{)}$ is n-embedded in G.}\par
\noindent\textbf{Corollary 4.11} \cite[Theorem C]{Guo1}. \textit{Let
$\mathfrak{F}$ be a saturated formation containing $\mathfrak{U}$
and $G$ a group with a normal subgroup E such that $G/E\in
\mathfrak{F}$. Suppose that for every non-cyclic Sylow subgroup P of
E, every maximal subgroup of P or every cyclic subgroup H of P of
prime order or order \textup{4} $\bigm{(}$if P is non-abelian
\textup{2}-group and $H\nsubseteq Z_{\infty}(G)$$\bigm{)}$ is
S-embedded in G.}\par \noindent\textbf{Corollary 4.12} \cite[Theorem
3.3]{Li2}. \textit{Let $\mathfrak{F}$ be a saturated formation
containing $\mathfrak{U}$. Then $G\in \mathfrak{F}$ if and only if
there exists a normal subgroup H of G such that $G/H\in
\mathfrak{F}$ and all maximal subgroups of any Sylow subgroup of
$F^*(H)$ are SS-quasinormal in G.}\par \noindent\textbf{Corollary
4.13} \cite[Theorem 3.7]{Li2}. \textit{Let $\mathfrak{F}$ be a
saturated formation containing $\mathfrak{U}$. Then $G\in
\mathfrak{F}$ if and only if G has a normal subgroup H such that
$G/H\in \mathfrak{F}$ and all cyclic subgroups of prime order or
order \textup{4} of $F^*(H)$ are SS-quasinormal in G.}\par
\noindent\textbf{Corollary 4.14} \cite[Theorem 2]{Zha}. \textit{Let
$\mathfrak{F}$ be a saturated formation containing $\mathfrak{U}$.
Then $G\in \mathfrak{F}$ if and only if there exists a soluble
normal subgroup H of G such that $G/H\in \mathfrak{F}$ and all
maximal subgroups of any Sylow subgroup of $F(H)$ are
S-semipermutable in G.}\par \noindent\textbf{Corollary 4.15}
\cite[Theorem 4]{Zha}. \textit{Let $\mathfrak{F}$ be a saturated
formation containing $\mathfrak{U}$. Then $G\in \mathfrak{F}$ if and
only if there exists a soluble normal subgroup H of G such that
$G/H\in \mathfrak{F}$ and all cyclic subgroups of prime order or
order \textup{4} of $F(H)$ are S-semipermutable in G.}\par
\noindent\textbf{Corollary 4.16} \cite[Theorem 4.3]{Wei}.
\textit{Let $\mathfrak{F}$ be a saturated formation containing
$\mathfrak{U}$. Suppose that $G$ is a group with a normal subgroup H
such that $G/H\in \mathfrak{F}$. If all maximal subgroups of any
Sylow subgroup of $F^*(H)$ are $c^*$-normal in G, then $G\in
\mathfrak{F}$.}\par \noindent\textbf{Corollary 4.17} \cite[Theorem
E]{Guo2}. \textit{Let $\mathfrak{F}$ be a saturated formation
containing $\mathfrak{U}$ and $G$ a group with a normal subgroup E
such that $G/E\in \mathfrak{F}$. If for every non-cyclic Sylow
subgroup P of $F^*(E)$, every maximal subgroup of P or every cyclic
subgroup H of P of prime order or order \textup{4} $\bigm{(}$if P is
non-abelian \textup{2}-group and $H\nsubseteq
Z_{\infty}(G)$$\bigm{)}$ is n-embedded in G, then $G\in
\mathfrak{F}$.}\par \noindent\textbf{Corollary 4.18} \cite[Theorem
D]{Guo1}. \textit{Let $\mathfrak{F}$ be a saturated formation
containing $\mathfrak{U}$ and $G$ a group with a normal subgroup E
such that $G/E\in \mathfrak{F}$. Suppose that for every non-cyclic
Sylow subgroup P of $F^*(E)$, every maximal subgroup of P or every
cyclic subgroup H of P of prime order or order \textup{4}
$\bigm{(}$if P is non-abelian \textup{2}-group and $H\nsubseteq
Z_{\infty}(G)$$\bigm{)}$ is S-embedded in G. Then $G\in
\mathfrak{F}$.}\par

\bibliographystyle{plain}
\bibliography{expbib}
\end{document}